\documentclass[12pt]{amsart}
\usepackage{geometry}                
\usepackage{graphicx}
\usepackage{amssymb}
\usepackage{epstopdf}
\usepackage{enumerate}
\newtheorem{df}{Definition}[section]

\newtheorem{thm}{Theorem}[section]
\newtheorem{prop}{Proposition}[section]

\newtheorem{ex}{Example}[section]
\newtheorem{remark}{Remark}[section]
\newtheorem{fact}{Fact}[section]
\newtheorem{cor}{Corollary}[section]

\title{Zeta functions of Ramanujan graphs and modular forms}
\author{Kennichi Sugiyama}

\begin{document}
\maketitle

\begin{center}
Department of Mathematics, Faculty of Science,\\ 
Rikkyo University, 3-34-1 Nishi-Ikebukuro, Toshima,\\
Tokyo 171-8501, Japan \\
e-mail address : kensugiyama@rikkyo.ac.jp
\end{center}

\begin{abstract} 
We will investigate the relationship between Ihara's zeta functions of Ramanujan graphs and Hasse-Weil's congruent zeta functions of modular curves. As an application we will describe the limit value of Hasse-Weil's congruent zeta functions in terms  of the corresponding Ramanujan graphs. Moreover we will show a congruence relation of the Fourier coefficients of a normalized Hecke eigenform of weight $2$.

AMS classification 2010: 05C30, 05C38, 05C50, 11G18, 11G20, 11M38, 11M99, 14G10.
\end{abstract}

\section{Introduction}
The aim of this report is to study the relationship between Ihara's zeta functions of connected Ramanujan graphs and Hasse-Weil zeta functions of modular curves on a certain prime level. \\

Shortly, a graph is a one dimensional simplicial complex. In this paper, we will only treat graphs whose edges are oriented in both directions. The information of a graph $G$ are encoded in the adjacency matrix $A$, which describes how edges and vertices are connected. Although graphs are geometric objects they are intimately related to number theory. In fact, let ${\mathcal P}(G)$ be the set of reduced and tail-less primitive closed paths of $G$ i.e. the set of closed path without backtracking and not around more than once. Two closed paths are defined to be {\it equivalent}, if one is obtained from the other by a cyclic shift of edges. Let $\frak{P}(G)$ denote the set of equivalence classes of ${\mathcal P}(G)$. 
The length of a closed path only depends on the equivalence class and we have a function
\[l : \frak{P}(G) \to {\mathbb Z}\]
and {\em Ihara's zeta function} of $G$ is defined as
\[Z(G;t)=\prod_{[c]\in {\frak P}(G)}\frac{1}{1-t^{l([c])}}\]
(See {\bf Section 2} for details of these materials and the facts concerning graphs). The function is originally defined by Ihara and has been studied by various mathematicians (we only refer to \cite{Hashimoto1990}, \cite{Hashimoto1991}, \cite{Ihara1965a}, \cite{Ihara1965b}, \cite{Ihara1966}, \cite{Hoffman}, \cite{Stark-Terras} but there are much more). One of the most remarkable properties of $Z(G;t)$ is that it is a rational function:
\begin{equation}Z(G;t)=\frac{(1-t^2)^{\chi(G)}}{{\rm det}[1-At+Qt^2]}.\end{equation}
Here $Q$ is a diagonal matrix whose entry at the vertex $x$ is $d(x)-1$ where $d(x)$ is the number of edges exiting from $x$, and $\chi(G)$ denotes the Euler characteristic of $G$. In this paper we assume that $G$ is connected and $d$-regular i.e. $d(x)=d$ for any vertex $x$. Then $d$ is an eigenvalue of $A$ of multiplicity one. Moreover, it is known that the modulus of an eigenvalue $\lambda$ of $A$ is less than or equal to $d$, and that $-d$ is an eigenvalue of $A$ if and only if $G$ is bipartite. If $|\lambda|\leq 2\sqrt{d-1}$ is satisfied for any eigenvalue $\lambda$ of $A$ other than $\pm d$, the graph is called {\it Ramanujan}.\\

 We recall that Hasse-Weil zeta functions have similar properties as  (1). Let $C$ be a smooth proper curve defined over a finite field ${\mathbb F}_q$ of characteristic $p$. The $q$-th power arithmetic Frobenius $\phi_q$ acts on the set of closed points $C(\overline{\mathbb F}_q):={\rm Hom}_{{\mathbb F}_q}({\rm Spec}(\overline{\mathbb F}_q), C)$ in the obvious way. Let $|C|$ be the orbit space. The degree of a closed point $x$ of $C$ is defined to be the extension degree of the residue field $k(x)$ over ${\mathbb F}_q$ and it factors through the map
\[{\rm deg}\,:\, |C| \to {\mathbb Z}.\]
The {\em Hasse-Weil zeta function} $W(C;t)$ of $C$ is defined as
\begin{equation}W(C;t)=\prod_{x\in|C|}\frac{1}{1-t^{{\rm deg}(x)}}.\end{equation}
{\it A priori} this is only a formal power series of $t$, but it turns out that it is a rational function by the Grothendieck-Lefschetz trace formula (\cite{SGA41/2}). In fact for a prime $l(\neq p)$ let $H^i_{et}(C_{{\overline{\mathbb F}}_q}, {\mathbb Z}_l)$ denote the $l$-adic \'{e}tale cohomology, a free ${\mathbb Z}_l$-module whose rank is the twice of the genus of $C$. Here $C_{{\overline{\mathbb F}}_q}$ is the base change $C$ over an algebraic clusure $\overline{\mathbb F}_q$. The $q$-th power geometric Frobenius $Fr_q$ acts on $H^1_{et}(C_{{\overline{\mathbb F}}_q}, {\mathbb Z}_l)$, and 
\begin{equation}W(C;t)=\frac{{\rm det}(1-\rho_l(Fr_q)t\,|\,H^1_{et}(C_{{\overline{\mathbb F}}_q}, {\mathbb Z}_l))}{(1-t)(1-qt)}.\end{equation}

Comparing (1) and (2), it is natural to expect that there should be relation between Ihara's zeta functions and Hasse-Weil's zeta functions. In fact, we will construct Ramanujan graphs whose Ihara's zeta functions are intimately related to Hasse-Weil's zeta functions of modular curves. This is pointed out by Ihara (\cite{Ihara1967}) and after that there are several ways to construct Ramanujan graphs (\cite{Li1992}, \cite{L-P-S}, \cite{Mestre}). Our construction is based on the way of Mestre and Oesterl\'{e}, which is sketched in \cite{Mestre}. Here is a summary of our construction. Let us fix a prime $N$ and take another prime $p$. The Hecke operator $T_p$ naturally acts on the free abelian group generated by supersingular elliptic curves over $\overline{\mathbb F}_N$. The representation matrix is called {\em a Brandt matrix} and will be denoted by $B(p)$. As we will see in the next section, taking adjacency matrix gives a bijective correspondence between the set of graphs and the set of matrices satisfying certain conditions ({\bf Proposition 2.1}). The Brandt matrix is a candidate of an adjacency matrix of our graph but unfortunately not in general since it does not satisfy a necessary condition.  However, if $N-1$ is divisible by $12$, $B(p)$ becomes an adjacency matrix of a graph $G_N(p)$ ({\bf Proposition 3.1}). In the remaining of this section, we assume that $N-1$ is divisible by $12$, and we set
\[n:=\frac{N-1}{12}.\]
Then $G_N(p)$ turns out to be a connected $(p+1)$-regular Ramanujan graph which is not bipartite ({\bf Theorem 3.1}). Moreover we will show that Ihara's zeta function of $G_{N}(p)$ and Hasse-Weil's zeta function of $X_{0}(N)_{{\mathbb F}_p}$ are {\em reciprocal}. 
\begin{thm}
\begin{enumerate}
\item 
\[W(X_{0}(N)_{{\mathbb F}_p},t)Z(G_N(p),t)=\frac{1}{(1-t)^2(1-pt)^2(1-t^2)^{\frac{n(p-1)}{2}}}.\]
\item 
\[\lim_{t\to 1}(t-1)W(X_{0}(N)_{{\mathbb F}_p},t)=\frac{n\tau(G_N(p))}{p-1},\]
where $\tau(G_N(p))$ is the complexity of $G_N(p)$ (see \S 2).
\end{enumerate}
\end{thm}
The assertion (2) may be compared to {\em the class number formula}, but we do not know how $\tau(G_N(p))$ can be interpreted as a $K$-theoretical object. Let $S_2(\Gamma_0(N))$ be the space of cusp forms of weight $2$ for the Hecke congruence subgroup $\Gamma_0(N)$ of ${\rm SL}_2({\mathbb Z})$. Then the dimension of $S_2(\Gamma_0(N))$ is $n-1$ and we can take normalized Hecke eigenforms $\{f_1, \cdots, f_{n-1}\}$ as its basis. Let
\[f_i=\sum_{n=1}^{\infty}a_n(f_i)q^n, \quad q=e^{2\pi iz}\quad({\rm Im}\,z>0)\]
be the Fourier expansion which satisfies $a_1(f_{i})=1$ by definition. By the classical Kirchhoff formula ({\bf Proposition 2.3}) which is used to derive (2) of {\bf Theorem 1.1}, and Eichler-Shimura's relation we will show the following congruence.

\begin{thm} Let $p(\neq N)$ be a prime such that $1+p$ is a multiple of $n$. Then the product of $p$-th coefficients $\mu_{N}(p):=\prod_{i=1}^{n-1}a_{p}(f_i)$ is divisible by $n$.
\end{thm}

We have listed the results of numerical experiments for $N=37,\,61,\,73$ in \S 5. {\bf Theorem 1.2} immediately yields the following corollary.
\begin{cor} Let $r$ be a prime divisor of $n$. Then there is a normalized Hecke eigenform $f=\sum_{n=1}^{\infty}a_n(f)q^n \in S_2(\Gamma_0(N))$ satisfying
\[|\{p\,\,\text{is a prime} : a_p(f)\equiv 0\, ({\rm mod}\,\frak{R}), \quad p\equiv -1\,({\rm mod}\,n)\}|=\infty\]
where $\frak{R}$ is a prime ideal of the integer ring of ${\mathbb Q}[\{a_n(f)\}_n]$ lying over $r$.
\end{cor}
{\bf Proof.} In fact we will show that at least one of $\{f_i\}_{i=1,\cdots,n-1}$ satisfies the assertion. Suppose that the claim were not true for all $\{f_i\}_{i=1,\cdots,n-1}$. Then except for finite primes $p$ with $p\equiv -1\,({\rm mod}\,n)$, 
$\prod_{i=1}^{n-1}a_{p}(f_i)$ is not divisible by $r$. This contradicts {\bf Theorem 1.2}. 
\begin{flushright}
$\Box$
\end{flushright}

{\bf Acknowledgements.} The author thanks Prof. Noro who has kindly informed us of the results of numerical experiments, and Prof. Geisser for his careful reading the manuscript. He also appreciates the referee's kind comments and valuable suggestions.

\section{The zeta function of a graph}
In this section we will recall the basic facts concerning Ihara's zeta functions of finite graphs and Ramanujan graphs, mostly following \cite{Bass} and \cite{Terras}.\\

A (finite) graph $G$ consists of a finite set of vertices $V(G)$ and a finite set of oriented edges $E(G)$, which satisfies the following property: there are {\it end point maps},
\[\partial_0, \quad \partial_1 : E(G) \to V(G),\]
and {\it an orientation resersal},
\[J:E(G) \to V(G),\quad J^2=\text{identity},\]
such that $\partial_i\circ J=\partial_{1-i}\,(i=0,1)$. The quotient $E(G)/J$ is called {\it the set of geometric edges} and is denoted by $GE(G)$. We regard an element of $GE(G)$ as an unoriented edge.  For $x\in V(G)$ we set
\[E_j(x)=\{e\in E(G)\, | \, \partial_j(e)=x\}, \quad j=0,1.\]
Thus $JE_j(x)=E_{1-j}(x)$. Intuitively $E_0(x)$ (resp. $E_1(x)$) is the set of edges starting from (resp. arriving at) $x$. The {\it degree} of $x$, $d(x)$, is defined by
\[d(x)=|E_0(x)|=|E_1(x)|.\]
$E(G)$ is naturally divided into two classes, {\it loops} and {\it  passes}. An edge $e\in E(G)$ is called {\it a loop} if $\partial_0(e)=
\partial_1(e)$ and is called {\it a pass} otherwise. Let $p(x)$ (resp. $l(x)$) be the number of passes (resp. the half of the number of loops, that is, forgetting the orientation, $l(x)$ is the number of {\em geometric} loops) starting from $x$.
Note that, because of the involution $J$, if we replace "starting" by "arriving" these number does not change. By definition, it is clear that
\[d(x)=2l(x)+p(x).\]
We set $q(x):=d(x)-1$. 
Let $C_{0}(G)$  be the free ${\mathbb Z}$-module generated by $V(G)$ with vertices as the natural basis. We define endomorphisms $Q$ and $A$ of $C_0(G)$ by
\[\quad Q(x)=q(x)x, \quad x\in V(G),\]
and 
\[A(x)=\sum_{e\in E(G), \partial_0(e)=x}\partial_1(e),\quad x\in V(G),\]
respectively. Note that because of the involution $J$, 
\[A(x)=\sum_{e\in E(G), \partial_1(e)=x}\partial_0(e).\]
The operator $A$ will be called {\it adjacency operator}.  We sometimes identify it with the representing matrix with respect to the basis $\{x\}_{x\in V(G)}$. Thus the $xy$-entry of $A_{xy}$ is the number of edges starting from $x$ and arriving at $y$. The orientation reversing involution $J$ implies
\[A_{xy}=A_{yx}. \]
Note that $A_{xx}=2l(x)$ and $p(x)=\sum_{y\neq x}A_{yx}$. If $d(x)=k$ for all $x\in V(G)$, $G$ is called $k$-{\it regular}.\\

Connecting distinct vertices $x$ and $y$ by $A_{xy}$-edges and drawing $\frac{1}{2}A_{xx}$-loops at $x$,  the adjacency matrix $A$ determines a $1$-dimensional unoriented simplicial complex. We call it {\it the geometric realization} of $G$, and denote it by $G$ again. We say that $G$ is connected if its the geometric realization is. The Euler characteristic $\chi(G)$ is equal to $|V(G)|-|GE(G)|$, hence if $G$ is connected, the fundamental group is a free group of rank $1-|V(G)|+|GE(G)|$. {\it A tree} is defined to be a graph which is connected and simply connected. A tree $T$ contained in $G$ satisfying $V(T)=V(G)$ is called {\it a spanning tree} of $G$.  Let $\tau(G)$ denote the number of spanning trees of $G$. We call $\tau(G)$ the {\it complexity} of $G$. For a later purpose, we summarize the relationship between a graph and its adjacency matrix.
\begin{prop} Let $A=(a_{ij})_{1\leq i,j \leq m}$ be an $m\times m$-matrix satisfying the following conditions.
\begin{enumerate}[(1)]
\item The entries $\{a_{ij}\}_{ij}$ are non-negative integers and satisfy
\[a_{ij}=a_{ji}, \quad \forall i \,\text{and}\, j.\]
\item $a_{ii}$ is even for every $i$.
\end{enumerate}
Then there is a unique graph $G$ whose adjacency matrix is $A$. Moreover, $G$ is $k$-regular if and only if one of the following conditions holds :
\begin{enumerate}[(a)]
\item
\[\sum_{i=1}^{m}a_{ij}=k,\quad \forall j\]
\item 
\[\sum_{j=1}^{m}a_{ij}=k,\quad \forall i.\]
\end{enumerate}
\end{prop}
 
\begin{ex} Let $G$ and $G^{\prime}$ be graphs with vertices $\{[1],[2],[3],[4]\}$ and whose shape are described by the following adjacency matrices $A$ and $A^{\prime}$, respectively. The $(i,j)$-entry is the number of geometric edges connecting vertices $[i]$ and $[j]$.
\begin{enumerate}[(I)]
\item
\[A=\left(
  \begin{array}{cccc}
    2 & 1 & 0 & 1 \\ 
    1 & 0 & 3 & 0 \\ 
    0 & 3 & 0 & 1 \\ 
    1 & 0 & 1 & 2 \\ 
  \end{array}
\right), \quad 
Q=\left(
  \begin{array}{cccc}
    3 & 0 & 0 & 0 \\ 
    0 & 3 & 0 & 0 \\ 
    0 & 0 & 3 & 0 \\ 
    0 & 0 & 0 & 3 \\ 
  \end{array}
\right)
\]
\item 
\[A^{\prime}=\left(
  \begin{array}{cccc}
    0 & 1 & 0 & 1 \\ 
    1 & 0 & 3 & 0 \\ 
    0 & 3 & 0 & 1 \\ 
    1 & 0 & 1 & 0 \\ 
  \end{array}
\right), \quad 
Q^{\prime}=\left(
  \begin{array}{cccc}
    1 & 0 & 0 & 0 \\ 
    0 & 3 & 0 & 0 \\ 
    0 & 0 & 3 & 0 \\ 
    0 & 0 & 0 & 1 \\ 
  \end{array}
\right)
\]
\end{enumerate}
 $G$ is a $4$-regular graph which has loops at the vertices $[1]$ and $[4]$. $G^{\prime}$ is obtained by deleting these loops from $G$. 
\end{ex}

In the following, a graph $G$ is always assumed to be {\it connected}. {\it A path of length} $m$ is a sequence $c=(e_1,\cdots,e_m)$ of edges such that $\partial_0(e_i)=\partial_1(e_{i-1})$ for all $1< i \leq m$ and the path is {\it reduced} if $e_i\neq J(e_{i-1})$ for all      $1< i \leq m$. The path is {\it closed} if $\partial_0(e_1)=\partial_1(e_{m})$, and the closed path has {\it no tail} if $e_m\neq J(e_1)$. A closed path of length one is nothing but a loop. Two closed paths are {\it equivalent} if one is obtained from the other by a cyclic shift of the edges. Let ${\frak C}(G)$ be the set of equivalence classes of reduced and tail-less closed paths of $G$. Since the length depends only on the equivalence class, the length function descends to the map;
\[l : {\frak C}(G) \to {\mathbb N},\quad l([c])=l(c),\]
where $[c]$ is the class determined by $c$. We define a reduced and tail-less closed path $C$ to be primitive if it is not obtained by going $r\,(\geq 2)$ times some another closed path. Let ${\frak P}(G)$ be the subset of ${\frak C}(G)$ consisting of the classes of primitive closed paths (which are  reduced and tail-less by definition). The {\em Ihara zeta function} of $G$ is defined to be
\[Z(G;t)=\prod_{[c]\in {\frak P}(G)}\frac{1}{1-t^{l([c])}}.\]
Although this is an infinite product, it is a rational function.
\begin{fact}(\cite{Bass},\cite{Hoffman},\cite{Ihara1966},\cite{Stark-Terras})
\[Z(G;t)=\frac{(1-t^2)^{\chi(G)}}{{\rm det}[1-At+Qt^2]}.\]
\end{fact}
\begin{prop}
Let $G$ be a $k$-regular graph with $m$ vertices. Then the Euler characteristic $\chi(G)$ satisfies
\[\chi(G)=\frac{m(2-k)}{2}.\]
\end{prop}
\begin{remark}
Note that the Euler characteristic does not depend on the number of loops.
\end{remark}
{\bf Proof.} Let $\{[1],\cdots,[m]\}$ be the vertices, and denote the number of distinct primitive loops passing through $[i]$ (resp. edges connecting distinct vertices $[i]$ and $[j]$) by $l(i)$ (resp. $r_{ij}$). Since $G$ is $k$-regular we know that
\[2l(i)+\sum_{j=1,j\neq i}^{n}r_{ij}=k,\quad 1\leq \forall i \leq m.\]
Take these sum and we find
\begin{equation}2\sum_{i=1}^{n}l(i)+\sum_{i,j, i\neq j}r_{ij}=mk.\end{equation}
Note that the number of edges $e$ is given by
\[e=\sum_{i=1}^{n}l(i)+\frac{1}{2}\sum_{i,j, i\neq j}r_{ij}.\]
Thus (4) implies 
\[e=\frac{1}{2}mk.\]
and 
\[\chi(G)=m-e=\frac{m(2-k)}{2}.\]
\begin{flushright}
$\Box$
\end{flushright}
\begin{cor}
Let $G$ be a $k$-regular graph with $m$ vertices. Then
\[Z(G;t)=\frac{(1-t^2)^{\frac{m(2-k)}{2}}}{{\rm det}[1-At+Qt^2]}.\]
\end{cor}

Let $E_{or}(G)\subset E(G)$ be a section of the natural projection $E(G) \to GE(G)$. In other word we choose an orientation on geometric edges and make the geometric realization into an oriented one dimensional simplicial complex. Let $C_1(G)$ be the free ${\mathbb Z}$-module generated by $E_{or}(G)$. Then the boundary map 
\[\partial : C_1(G) \to C_0(G)\]
is naturally defined. Let $\partial^{t}$ be the dual of $\partial$. The {\it Laplacian} $\Delta$ of $G$ is defined to be $\Delta=\partial \partial^{t}$. It is known (and easy to check) that (\cite{Terras}, \cite{Hoffman}),
\begin{equation}\Delta=1-A+Q.\end{equation}
Since we have assumed that $G$ is connected,  $0$ is an eigenvalue of $\Delta$ with multiplicity one. \cite{Biggs}{\bf Theorem 6.3} implies the following result.
\begin{prop}(The Kirchhoff law) Let $m$ be the number of vertices and $\{0, \lambda_1,\cdots, \lambda_{m-1}\}$ the eigenvalues of $\Delta$. Then the complexity is given as
\[\tau(G)=\frac{\lambda_1\cdots \lambda_{m-1}}{m}.\]
\end{prop}
Let $G^{\prime}$ be the graph obtained by deleting all loops of $G$. By the equation (5) we see that the Laplacian of $G$ and $G^{\prime}$ are same, and {\bf Proposition 2.1} implies
\[\tau(G)=\tau(G^{\prime}).\] 
One may also observe this fact by inspection.

\begin{ex} Let $G$ and $G^{\prime}$ be the graphs considered in {\bf Example 2.1}. The equation (1) shows that the corresponding Laplacians $\Delta$ and $\Delta^{\prime}$ are equal as expected;
\[\Delta=\Delta^{\prime}=
\left(
  \begin{array}{cccc}
    2 & -1 & 0 & -1 \\ 
    -1 & 4 & -3 & 0 \\ 
    0 & -3 & 4 & -1 \\ 
    -1 & 0 & -1 & 2 \\ 
  \end{array}
\right).\]
The eigenvalues of this matrix are $\{0,2,5+\sqrt{5}, 5-\sqrt{5}\}$, and Kirchhoff law tells us
\[\tau(G)=\tau(G^{\prime})=\frac{2(5+\sqrt{5})(5-\sqrt{5})}{4}=10.\] 
Drawing a picture, one can immediately verify this.
\end{ex}
Now let $G$ be a connected $k$-regular graph. As we have seen, since $0$ is an eigenvalue of $\Delta$ with multiplicity one, (5) shows that $k$ is an eigenvalue of $A$ with multiplicity one. Because of semi-positivity, $\Delta$ we find that 
\[|\lambda|\leq k\quad \text{for any eigenvalue $\lambda$ of $A$}\]
and that $-k$ is an eigenvalue of $A$ if and only if $G$ is bipartite (\cite{Terras}, {\bf Chapter 3}). Here $G$ is called {\em bipartite} if the set of vertices $V(G)$ can be divided into disjoint subset $V_0$ and $V_1$ such that every edge connects points in $V_0$ and $V_1$, namely there is no edge whose end points are simultaneously contained in $V_{0}$ or $V_{1}$. 

\begin{df} Let $G$ be a $k$-regular graph. We say that it is {\rm Ramanujan}, if all eigenvalues $\lambda$ of $A$ with $|\lambda|\neq k$ satisfy
\[|\lambda|\leq 2\sqrt{k-1}.\]
\end{df}
See \cite{Li1992}, \cite{Murty} and \cite{Valette} for detailed expositions of Ramanujan graphs.

\section{Ramanujan graphs associated to modular curves}
\subsection{The Brandt matrix}
In this subsection  we will recall the theory of Brandt matrices after \cite{Gross}. Let $N$ be a prime, and let $B$ be the quaternion algebra over ${\mathbb Q}$ ramified at two places $N$ and $\infty$. Let $R$ be a fixed maximal order in $B$ and $\{I_1,\cdots,I_n\}$ be the set of left $R$-ideals representing the distinct ideal classes. We call $n$ {\it the class number} of $B$. We choose $I_1=R$. For $1\leq i \leq n$, $R_i$ denotes the right order of $I_i$, and let $w_i$ the order of $R_i^{\times}/\{\pm 1\}$. 
The product 
\begin{equation}W=\prod_{i=1}^{n}w_{i}\end{equation}
is independent of the choice of $R$ and is equal to the exact denominator of $\frac{N-1}{12}$ (\cite{Gross}, p.117). Eichler's mass formula states that
\[\sum_{i=1}^{n}\frac{1}{w_i}=\frac{N-1}{12}.\]
Let ${\mathbb F}$ be an algebraic closure of ${\mathbb F}_N$. There are $n$ distinct isomorphism classes $\{E_1,\cdots,E_n\}$ of supersingular elliptic curves over ${\mathbb F}$ such that ${\rm End}(E_i)\simeq R_i$. Let $p$ be a prime distinct from $N$, and let ${\rm Hom}(E_i,\,E_j)(p)$ denote the set of homomorphisms from $E_i$ to $E_j$ of degree $p$. The $(i,j)$-entry of the Brandt matrix $B(p)$ is defined to be
\begin{equation}b_{ij}=\frac{1}{2w_j}|{\rm Hom}(E_i,\,E_j)(p)|.\end{equation}
Since ${\rm Hom}(E_i,\,E_j)(p)$ has a faithful action of $R_{j}^{\times}$ from the right, $b_{ij}$ is a non-negative integer. In fact $b_{ij}$ equals to the number of subgroup $C$ of order $p$ in $E_i$ such that $E_{i}/C\simeq E_j$ (\cite{Gross} {\bf Proposition 2.3}). \\

Now we assume that $N-1$ is divisible $12$. Since $\frac{N-1}{12}$ is an integer $W=\prod_{i=1}^{n}w_{i}=1$ and $w_i=1$ for all $i$. Hence by Eichler's mass formula 
\begin{equation}
n=\frac{N-1}{12}.
\end{equation}
Although the following statements are included in \cite{Pizer} {\bf Proposition 4.6}, we will give a proof for the sake of convenience.
\begin{prop} Let $N$ be a prime such that $N-1$ is divisible by $12$. Then the Brandt matrix $B(p)=(b_{ij})_{1\leq i,j \leq n}\, (p\neq N)$ satisfies the following.
\begin{enumerate}[(1)]
\item Every entry is a non-negative integer and $B(p)$ is symmetric;
\[b_{ij}=b_{ji}.\]
\item
The diagonal entires $\{b_{ii}\}_{i}$ are even for all $i$.
\item For any $i=1,\cdots, n$,
\[\sum_{j=1}^{n}b_{ij}=p+1.\]
\end{enumerate}
\end{prop}
{\bf Proof.} Taking the dual homomorphisms, we have a bijective correspondence
\[I : {\rm Hom}(E_i,\,E_j)(p) \to {\rm Hom}(E_j,\,E_i)(p)\]
defined by $I(\phi)=\check{\phi}$.
Since $w_i=w_j=1$ the claim (1) is clear from the definition. In order to show the claim (2), it is sufficient to show that the action of $I$ on ${\rm End}(E_i)(p)/{\pm 1}$ has no fixed point. Let $\phi$ be an element of ${\rm End}(E_i)(p)/{\pm 1}$. Then the kernels of $\phi$ and $\check{\phi}$ generate the subgroup of $p$-torsion points $E_i[p]$ :
\[{\rm Ker}(\phi)+{\rm Ker}(\check{\phi})=E_i[p]\simeq {\mathbb F}_p\oplus {\mathbb F}_p.\]
If $\phi\in {\rm End}(E_i)(p)/{\pm 1}$ were a fixed point of $I$, then ${\rm Ker}(\phi)={\rm Ker}(\check{\phi})$ and 
\[{\mathbb F}_p\simeq {\rm Ker}(\phi)=E_i[p]\simeq {\mathbb F}_p\oplus {\mathbb F}_p,\]
which is absurd. The claim (3) follows from the following observation : by definition $\sum_{j=1}^{n}b_{ij}$ is equal to the number of cyclic subgroups of $E_i$ of order $p$, which is $p+1$.
\begin{flushright}
$\Box$
\end{flushright}

By {\bf Proposition 2.1} $B(p)$ determines a $(p+1)$-regular graph, which will be denoted by $G_{N}(p)$. Here is a remark. In our proof of {\bf Proposition 3.1}, the fact that $R_i^{\times}=\{\pm 1\}\,(\forall i)$ plays a crucial role. Indeed if $N$ does not satisfy this assumption, then $\frac{N-1}{12}$ is not  an integer and there is some $i$ with $w_i>1$. Therefore the assumption is necessary.

\subsection{Modular curves and Ramanujan graphs}
In this section we will relate $G_{N}(p)$ with the space of modular forms. Let $N$ be a prime and $S_2(\Gamma_0(N))$ the space of cusp forms of weight $2$ for the Hecke's congruence subgroup
\[\Gamma_{0}(N):=\{\left(
  \begin{array}{cc}
    a & b \\ 
    c & d \\ 
  \end{array}\right) \in {\rm SL}_2({\mathbb Z}) : \quad c \equiv 0\, ({\rm mod}\,N)\}.
\]
Let $Y_{0}(N)$ be the modular curve which parametrizes isomorphism classes of a pair ${\bf E}=(E,\Gamma_N)$ of an elliptic curve $E$ and a cyclic subgroup $\Gamma_N$ of order $N$ of $E$. It is a smooth curve defined over ${\mathbb Q}$, and the set of ${\mathbb C}$-valued points is the quotient of the upper half plane by 
$\Gamma_0(N)$. The compactification $X_0(N)$ of $Y_0(N)$ has a canonical model over ${\mathbb Z}$ which has been studied by \cite{Deligne-Rapoport} and \cite{Katz-Mazur} in detail. Then $S_2(\Gamma_0(N))$ is identified with the space of holomorphic $1$-forms $H^0(X_0(N),\Omega)$, and in particular  with the tangent space ${\rm Tan}J_0(N)$ at the origin of the Jacobian variety $J_0(N)$ of  $X_0(N)$.\\

For a prime $p$ different from $N$, $X_0(N)$ furnishes the $p$-th Hecke operator defined by
 \begin{equation}T_p(E,\Gamma_N):=\sum_{C}(E/C, (\Gamma_N + C)/C),\end{equation}
where $C$ runs through all cyclic subgroups of $E$ of order $p$. By functoriality, $T_p$ acts on $J_0(N)$ and in turn on ${\rm Tan}J_0(N)$, and the action coincides with the usual action on $S_2(\Gamma_0(N))$ (see  \cite{Shimura1971}) under identification. We define the Hecke algebra as ${\mathbb T}:={\mathbb Q}[\{T_p\}_{p\neq N}\}]$, which is a commutative subring of ${\rm End}J_0(N)$. 
Let ${\mathcal J}_0(N)$ be the Neron model of $J_0(N)$ over ${\mathbb Z}$. It is known that  the connected components of the reduction of ${\mathcal J}_0(N)$ at $N$ are tori.   Following \cite{Ribet1990}, we will describe its character group. $X_0(N)_{{\mathbb F}_N}$ has two irreducible components $C_F$ and $C_V$, which are isomorphic to the projective line ${\mathbb P}^{1}=X_0(1)$. Over $C_F$ (resp. $C_V$), $\Gamma_N$ is the kernel of the Frobenius $F$ (resp. the Verschiebung $V$), and they intersect at supersingular points $\Sigma_N=\{E_1,\cdots,E_n\}$ as ordinary double points. Consider the homomorphism
\[\partial : \oplus_{E_i\in \Sigma_N}{\mathbb Z}E_i \to {\mathbb Z}C_F\oplus {\mathbb Z}C_V, \quad \partial(E_i)=C_F-C_V.\]
The image of $\partial$ is a free abelian group of rank one generated by $\delta:=C_F-C_V$. 
$X$ being the kernel of $\partial$, we have the exact sequence
\begin{equation}
0 \to X \to \oplus_{E_i\in \Sigma_N}{\mathbb Z}E_i \stackrel{\partial}\to {\mathbb Z}\delta \to 0.
\end{equation}
It is straightforward to check that 
\[X=\{\sum_{E_i\in \Sigma_N}a_i\cdot E_i\,|\,a_i\in{\mathbb Z},\,\sum_{i=1}^{n}a_i=0\},\]
and by \cite{Ribet1990}{\bf Proposition 3.1}, $X$ is the character group of the connected component of the reduction ${\mathcal J}_0(N)_{{\mathbb F}_N}$ at $N$. In particular we see that
\[{\rm dim}J_0(N)={\rm dim}S_2(\Gamma_0(N))={\rm dim}X\otimes{\mathbb Q}=n-1.\]
Let $p$ be a prime with $p\neq N$. Then $T_p$ operates on $\oplus_{E_i\in \Sigma_N}{\mathbb Z}E_i$ by (9). Note that  $\Gamma_N=0$, since $E_i$ is supersingular and we may write (9) as
 \begin{equation}T_p(E_i):=\sum_{C}E_i/C.\end{equation}
A simple computation shows that
\begin{equation}T_p(C_F)=(p+1)C_F,\quad T_p(C_V)=(p+1)C_V,\quad T_{p}(\delta)=(p+1)\delta.\end{equation}
and $X$ is preserved by $T_p$.  Remember that $X$ is the character group of ${\mathcal J}_0(N)_{{\mathbb F}_N}$ and this action of $T_p$ on $X$ is nothing but the action which is induced by $T_p$ on ${\mathcal J}_0(N)_{{\mathbb F}_N}$. Here is a relationship between $T_p$ and the Brandt matrix. By \cite{Gross}{\bf Proposition 4.4}
\begin{equation}T_pE_i=\sum_{j=1}^{n}b_{ij}E_j,\end{equation}
hence $B(p)$ is the representation matrix of $T_p$. It seems that the following statement is fairly well-known, but since we could not find out a suitable reference, we will give a proof.
\begin{prop}
$X\otimes{\mathbb C}$ and $S_2(\Gamma_0(N))$ are isomorphic as ${\mathbb T}$-modules.
\end{prop}
{\bf Proof.} Since the action of the Hecke operators $\{T_p\}_{p\neq N}$ on $S_2(\Gamma_0(N))$ is self-adjoint with respect to the Petersson inner product 
and is commutative, $\{T_p\}_{p\neq N}$ are simultaneously diagonalizable. Namely, there is a set of normalized Hecke eigenforms $\{f_1,\cdots, f_{n-1}\}$ which is a basis of $S_2(\Gamma_0(N))$, and for each pair $i,\,p$ there exists $\sigma_{i}(T_p)\in {\mathbb C}$ such that
\[T_p(f_i)=\sigma_{i}(T_p)f_i.\]
Thus $\sigma_i$ yields an algebraic homomorphism from ${\mathbb T}\otimes{\mathbb C}$ to ${\mathbb C}$, and by {\it the multiplicity one theorem} (\cite{Atkin-Lehner}), $\{\sigma_i\}_{i=1,\cdots, n-1}$  are mutually distinct. The perfect pairing
\[(\,,\,): {\mathbb T}\otimes_{\mathbb C} S_2(\Gamma_0(N)) \to {\mathbb C}, \quad (T,f)=a_1(Tf),\]
implies that the dimension of ${\mathbb T}\otimes{\mathbb C}$ is $n-1$, and we see that ${\rm Spec}({\mathbb T}\otimes{\mathbb C})({\mathbb C})=\{\sigma_1,\cdots,\sigma_{n-1}\}$. Thus there is an isomorphism of ${\mathbb C}$-algebras
\[{\mathbb T}\otimes{\mathbb C}\simeq {\mathbb C}^{n-1}\]
such that $\sigma_i$ is the projection to the $i$-th factor. 
 Now let us turn our attention to the character group $X$ of the connected component of ${\mathcal J}_0(N)_{{\mathbb F}_N}$. Define a positive definite inner product on $X\otimes{\mathbb R}$ to be
\[\langle E_i, E_j \rangle=w_i\delta_{ij},\]
where $\delta_{ij}$ is Kronecker's delta and the action of ${\mathbb T}$ is self-adjoint with respect to this (\cite{Gross} {\bf Proposition 4.6}).
Therefore $X\otimes{\mathbb C}$ has the spectral decomposition 
\[X\otimes{\mathbb C}=\oplus_{\sigma\in S}V_{\sigma},\]
where $S$ is a subset of ${\rm Spec}({\mathbb T}\otimes{\mathbb C})({\mathbb C})$ and $V_{\sigma}$ is one dimensional vector space over ${\mathbb C}$ which has a ${\mathbb T}$-action such that
\[T(x)=\sigma(T)x,\quad T\in{\mathbb T},\,x\in V_{\sigma}.\]
Since the action of ${\mathbb T}$ on $X\otimes {\mathbb C}$ is faithful (\cite{Ribet1990}{\bf Theorem 3.10}), every element of $\{\sigma_1,\cdots,\sigma_{n-1}\}$ appears in $S$. Indeed otherwise, the action factors through a quotient ${\mathbb T}/\frak{A}$ ($\frak{A}$ is a non-zero ideal), which contradicts faithfulness.  This observation and the equation ${\rm dim}X\otimes{\mathbb C}={\rm dim}S_2(\Gamma_0(N))=n-1$ imply that
\begin{equation}X\otimes{\mathbb C}\simeq S_2(\Gamma_0(N))\simeq \oplus_
{i=1}^{n-1}V_{\sigma_i}\end{equation}
as ${\mathbb T}$-modules.  

\begin{flushright}
$\Box$
\end{flushright}

Hence {\bf Proposition 3.2} together with (10) and (12) implies that 
\begin{eqnarray*}
{\rm det}[1-T_pt+pt^2\,|\,\oplus_{i}{\mathbb C}E_i] &=&  {\rm det}[1-T_pt+pt^2\,|\,X\otimes{\mathbb C}]\cdot {\rm det}[1-T_pt+pt^2\,|\,{\mathbb C}\delta]\\
&=& {\rm det}[1-T_pt+pt^2\,|\,S_2(\Gamma_0(N))](1-t)(1-pt).\\
\end{eqnarray*}
On the other hand, it follows from (13) that
\[{\rm det}[1-T_pt+pt^2\,|\,\oplus_{i}{\mathbb C}E_i] ={\rm det}[1-B(p)t+pt^2]. \]
Therefore,
\begin{equation}
{\rm det}[1-T_pt+pt^2\,|\,S_2(\Gamma_0(N))](1-t)(1-pt)={\rm det}[1-B(p)t+pt^2].
\end{equation}
It follows from {\bf Corollary 2.1} that
\[Z(G_N(p),t)=\frac{1}{{\rm det}[1-B(p)t+pt^2](1-t^2)^{\frac{n(p-1)}{2}}}.\]
Let $p$ be a prime not equal to $N$ such that $X_0(N)$ has a good reduction at $p$. Let us fix an arbitrary prime $l$ distinct from $p$.  Using Eichler-Shimura relation (\cite{Deligne}\cite{Shimura}), we see that the characteristic polynomial of the geometric Frobenius $Fr_p$ is computed by the Brandt matrix,
\begin{equation}
{\rm det}[1-B(p)t+pt^2]={\rm det}(1-p_l(Fr_p)t\,|\,H^{1}_{et}(X_0(N)_{\overline{{\mathbb F}}_p}, {\mathbb Z}_l)(1-t)(1-pt).
\end{equation}
Now we will relate these observations to Ihara's zeta function of a graph. In the following, let $N$ be a prime such that $N-1$ is divisible by $12$. As we have shown in {\bf Proposition 3.1}, the Brandt matrix $B(p)$ is the adjacency matrix of the $(p+1)$-regular graph $G_N(p)$. 

\begin{thm} Let $N$ be a prime such that $N-1$ is divisible by $12$ and $p$ a prime distinct from $N$. 
\begin{enumerate}[(1)]
\item
$G_N(p)$ is a connected $(p+1)$-regular Ramanujan graph which is not bipartite.
\item 
\[W(X_{0}(N)_{{\mathbb F}_p},t)Z(G_N(p),t)=\frac{1}{(1-t)^2(1-pt)^2(1-t^2)^{\frac{n(p-1)}{2}}}.\]

\item
\[{\rm det}(1-p_l(Fr_p)\,|\,H^{1}_{et}(X_0(N)_{\overline{{\mathbb F}}_p}, {\mathbb Z}_l))=\prod_{i=1}^{n-1}(1+p-a_{p}(f_i))=n\tau(G_N(p)).\]
\end{enumerate}
\end{thm}
{\bf Proof.}   
\begin{enumerate}
\item 
The claim (1) is an immediate consequence of (16) and the Weil conjectures for a curve. In fact, by (13) we know that $B(p)$ is the representation matrix of $T_p$, and $p+1$ is an eigenvalue of $B(p)$ by (12). By (9), the other eigenvalues $\{\lambda_1,\cdots, \lambda_{n-1}\}$ agree with the eigenvalues of the restriction of $T_p$ on $X\otimes{\mathbb C}$. Now taking into account of the Weil conjectures for a smooth curve, {\bf Proposition 3.2} and the Eichler-Shimura relation imply that the moduli of $\{\lambda_1,\cdots, \lambda_{n-1}\}$ are less than or equal to $2\sqrt{p}$. Thus we see that $G_N(p)$ is a Ramanujan graph which is not bipartite since $-(p+1)$ is not an eigenvalue of $B(p)$. Moreover, by (5) we see that the Laplacian $\Delta$ of $G_N(p)$ has zero eigenvalue with multiplicity one. Hence $G_N(p)$ is connected.
\item By (3) we know that 
\[W(X_0(N)_{{\mathbb F}_p};t)=\frac{{\rm det}(1-\rho_l(Fr_p)t\,|\,H^1_{et}(X_0(N)_{{\overline{\mathbb F}}_q}, {\mathbb Z}_l))}{(1-t)(1-pt)}.\]
and the claim (2) is an immediate consequence of {\bf Corollary 2.1} and (16).
 \item
 The left $"="$ of the claim (3) is a consequence of the Eichler-Shimura congruence relation. Indeed, equations (15) and (16) yield
\[{\rm det}(1-p_l(Fr_p)t\,|\,H^{1}_{et}(X_0(N)_{\overline{{\mathbb F}}_p}, {\mathbb Z}_l))=\prod_{i=1}^{n-1}(1-a_p(f_i)t+pt^2).\]
Setting $t=1$, the left equation is obtained. Next we will discuss the right $"="$ of the claim (3).  Remember that $B(p)$ is the representation matrix of $T_p$ and we have the exact sequence (10). Then the equations (12) and {\bf Proposition 3.2} show that the eigenvalues of $B(p)$ are 
\[\{p+1, a_p(f_1),\cdots, a_p(f_{n-1})\}.\]
Moreover, as we have seen, $B(p)$ is the adjacency matrix of a $(p+1)$-regular graph $G_N(p)$. Hence by (5) the eigenvalues of the Laplacian $\Delta$ of $G_N(p)$, which is $1+p-B(p)$, are 
 \[\{0, p+1-a_p(f_1),\cdots, p+1-a_p(f_{n-1})\}.\] 
Now the right equation follows from {\bf Proposition 2.3}.
\end{enumerate}
\begin{flushright}
$\Box$
\end{flushright}
{\bf Theorem 1.1} and {\bf Theorem 1.2} immediately follow from  {\bf Theorem 3.1}. Since we know that 
\[|a_p(f_i)| \leq 2\sqrt{p},\quad \forall i\]
by Weil's conjecture, {\bf Theorem 3.1 (3)} yields the following estimate.
\begin{cor}
\[\frac{(\sqrt{p}-1)^{2(n-1)}}{n} \leq \tau(G_N(p)) \leq \frac{(\sqrt{p}+1)^{2(n-1)}}{n}.\]
\end{cor}
\section{Numerical tables}
\begin{enumerate}[I]
\item Take $N=37$. Then $n=3$, and the dimension of $J_0(37)$ is two. Hence there are two cuspidal Hecke eigenform $f_{37,a}$ and $f_{37,b}$. Let $p$ be a prime such that $p+1$ is a multiple of $3$. {\bf Theorem 1.2} says that, for such a prime $p$, $\mu_{37}(p)$ is a multiple of $3$. Here are some calculations.

\begin{center}
\begin{tabular}{|c|c|c|c|c|c|c|c|} \hline
p & 5 & 11 & 17 & 23 & 29 & 41 & 47 \\ \hline\hline
$a_p(f_{37,a})$ &-2 &-5 &0 & 2&6 & -9&-9 \\ \hline
$a_p(f_{37,b})$ &0 &3 &6 & 6& -6& -9& 3\\ \hline\hline
$\mu_{37}(p)$ &0 &-15 &0 & 12& 36& 81& -27\\ \hline
\end{tabular}
\end{center}
\vspace{10mm}
\item 
Take $N=61$. Then $n=5$ and the dimension of $J_0(61)$ is $4$. Hence there are four cuspidal Hecke eigenform $f_{61a}$ and $\{f_{61b, (i)}\}_{i=1,2,3}$. Here $f_{61,a}$ is defined over ${\mathbb Q}$ and $\{f_{61b, (i)}\}_{i=1,2,3}$ are defined over $K$, where $K$ is the decomposition field of $x^3-x^2-3x+1=0$. Let $p$ be a prime such that $p+1$ is a multiple of $5$. By {\bf Theorem 1.2}, for such a prime $p$, $\mu_{61}(p)$ is a multiple of $5$. 

\begin{center}
\begin{tabular}{|c|c|c|c|c|c|} \hline
p & 19 & 29 & 59 & 79 & 89  \\ \hline\hline
$a_p(f_{61a})$ &4 &-6 & 9 & 3&-4  \\ \hline
$a_p(f_{61b, (i)})$ &$3\gamma_{i}-7$ & $-\gamma_{i}^2+2\gamma_{i}+3$ &$-\gamma_i^2-3\gamma_i+13$ & $-4\gamma_i^2-\gamma_i+14$ & $4\gamma_i^2-2\gamma_i-10$  \\ \hline\hline
$\mu_{37}(p)$ & 80 &120 &2925 & -1875& 320 \\ \hline
\end{tabular}
\end{center}
Here $\{\gamma_1, \gamma_2, \gamma_3\}$ are the distinct solutions of $x^3-x^2-3x+1=0$.
\vspace{10mm}
\item Take $N=73$. Then $n=6$ and the dimension of $J_0(73)$ is $5$. Hence there are five cuspidal Hecke eigenform $f_{73,a}$, $f_{73,b}^{\pm}$ and $f_{73,c}^{\pm}$. Although $f_{73,a}$ is defined over ${\mathbb Q}$, $f_{73,b}^{\pm}$ and $f_{73,c}^{\pm}$ are defined over ${\mathbb Q}(\sqrt{5})$ and ${\mathbb Q}(\sqrt{13})$, respectively. Here ${\pm}$ denotes the conjugate over ${\mathbb Q}$. Let $p$ be a prime such that $p+1$ is a multiple of $6$. {\bf Theorem 1.2} predicts that, for such a prime $p$, $\mu_{73}(p)$ is a multiple of $6$. 

\begin{center}
\begin{tabular}{|c|c|c|c|c|c|c|c|c|} \hline
p & 5 & 11 & 17 & 23 & 29 & 41 & 47 & 53 \\ \hline\hline
$a_p(f_{73,a})$ &2 &-2 &2 & 4& 2&6 & 6&10 \\ \hline
$a_p(f^{\pm}_{73,b})$ &$\alpha$ &$-\alpha-3$ &$-6\alpha-9$ &$\alpha-6$ & $-4\alpha-3$ &   $4\alpha+6$ & -$4\alpha-9$ & $8\alpha+15$ \\ \hline
$a_p(f^{\pm}_{73,c})$ & $-\beta$ & $\beta+3$ & $2\beta-3$ &$\beta+6$ &  $-4\beta+3$ & $-6$ & $9$ & $4\beta-3$ \\ \hline\hline
$\mu_{73}(p)$ &-6 &-18 & 810 &8580 & 1122& 720 & 396& 36210 \\ \hline
\end{tabular}
\end{center}
In the table $\alpha$ and $\beta$ are the solutions of
\[\alpha^2+3\alpha+1=0,\quad \beta^2-\beta-3=0.\]
\end{enumerate}


\end{document}